\documentclass[%
  onecolumn
]{mpi2015-cscpreprint}

\usepackage[american]{babel}
\usepackage{float}
\usepackage{cancel}
\usepackage{comment}
\usepackage{algpseudocode}
\usepackage{algorithm}
\usepackage{graphicx}
\usepackage{bm}
\usepackage{subcaption}
\usepackage{amssymb}
\usepackage{amsthm}
\theoremstyle{plain}
\newtheorem*{theorem*}{Theorem}

\usepackage[hardware]{mymacros}

\usepackage[normalem]{ulem}

\begin{document}
  

\title{On the representation of energy-preserving quadratic operators with application to Operator Inference}
  
\author[$1$,$\ast$]{Leonidas Gkimisis}
\author[$1$]{Igor Pontes Duff}
\author[$1$]{Pawan Goyal}
\author[$1$]{\newline Peter Benner}

\affil[$1$]{Max Planck Institute for Dynamics of Complex Technical Systems, Sandtorstra{\ss}e 1, 39106 Magdeburg, Germany,
Computational Methods in Systems and Control Theory (CSC).}
\affil[*]{Corresponding author: gkimisis@mpi-magdeburg.mpg.de}
  
\abstract{%
In this work, we investigate a skew-symmetric parameterization for energy-preserving quadratic operators. Earlier, [Goyal et al., 2023] proposed this parameterization to enforce energy-preservation for quadratic terms in the context of dynamical system data-driven inference. We here prove that every energy-preserving quadratic term can be equivalently formulated using a parameterization of the corresponding operator via skew-symmetric matrix blocks. Based on this main finding, we develop an algorithm to compute an equivalent quadratic operator with skew-symmetric sub-matrices, given an arbitrary energy-preserving operator. Consequently, we employ the skew-symmetric sub-matrix representation in the framework of non-intrusive reduced-order modeling (ROM) via Operator Inference (OpInf) for systems with an energy-preserving nonlinearity. To this end, we propose a sequential, linear least-squares (LS) problems formulation for the inference task, to ensure energy-preservation of the data-driven quadratic operator. The potential of this approach is indicated by the numerical results for a 2D Burgers' equation benchmark, compared to classical OpInf. The inferred system dynamics are accurate, while the corresponding operators are faithful to the underlying physical properties of the system.
}

\novelty{Skew-symmetric representation for any energy-preserving quadratic nonlinearity with application to non-intrusive model reduction}

\maketitle

\section{Introduction}%
\label{sec:intro}

Dynamical systems with quadratic nonlinearities are of specific interest to a wide range of engineering disciplines. A large class of autonomous dynamical systems can be lifted via a nonlinear transformation to a quadratic system \cite{Gu2011, Kramer2019, Benner2015}. Hence, examining the properties of quadratic system dynamics is especially valuable for control, simulation, and data-driven inference \cite{Kramer2021, Peherstorfer2016, Benner2015, Gosea2018, karachalios2022datadriven}. Here, we focus on systems with an energy-preserving quadratic nonlinearity. Intuitively, this means that the net contribution of the quadratic term to the system kinetic energy is zero. Thus, the quadratic term affects only how the kinetic energy is distributed to the system degrees of freedom (DOFs).

Although the above specification might seem too restrictive, it is satisfied by many nonlinear systems. Among others, a wide class of incompressible fluid dynamics problems falls under the category of dynamical systems with energy-preserving quadratic nonlinearity \cite{Schlegel2015, Lorenz1963, Holmes2012}. This property has driven stability analyses of incompressible fluid systems  \cite{Kalur2021, Schlegel2015} and has motivated physics-informed nonintrusive modeling methods \cite{Kaptanoglu2021, goyal2023guaranteed}. Considering the latter, \cite{Kaptanoglu2021} promoted the inference of globally stable quadratic models in the context of sparse identification \cite{Brunton2016}. Soft constraints on physics-based properties of the linear and quadratic operators were introduced to the objective function of a sparse, relaxed regularized regression. From a different perspective, the authors in \cite{goyal2023guaranteed} addressed the data-driven inference of both locally and globally stable systems by construction. This was achieved by encoding sufficient properties for the linear operator stability and the energy preservation of the quadratic operator to a gradient descent optimization task. For the quadratic operator, energy preservation was enforced by introducing a parameterization with skew-symmetric sub-matrices. 

Focusing on the quadratic term, the current work proves that any energy-preserving quadratic operator can be equivalently represented using the parameterization proposed in \cite{goyal2023guaranteed}. Moreover, we also showcase other possible representations for the energy-preserving nonlinearity. Based on these findings, we provide an algorithm that transforms any arbitrary energy-preserving operator $\mathbf{H}$ to an equivalent quadratic operator $\tilde{\mathbf{H}}$ with skew-symmetric sub-matrices. Finally, we employ this representation in the framework of physics-based, non-intrusive model reduction via OpInf. We propose a closed-form formulation for sequentially solving the involved linear LS problems to enforce energy preservation of the quadratic operator. We showcase the inference of an energy-preserving quadratic operator using the proven representation and the sequential LS approach for a 2D Burgers' test case. A comparison with a standard, unconstrained, sequential OpInf formulation illustrates that exploiting this physics-based property for systems with an energy-preserving quadratic term can result to qualitatively and quantitatively improved inferred models, for negligible added computational cost and coding effort.

\section{Representation of energy-preserving quadratic matrices}
\label{sec:fact}

We consider a quadratic dynamical system of the form

\begin{equation}
\label{quad}
\dot{\bx}(t) =\mathbf{A} \bx(t) + \mathbf{H} \left( \bx(t) \otimes \bx(t) \right) + \mathbf{B} \bu(t) + \mathbf{c},
\end{equation}

\noindent with state $\bx(t) \in \mathbb{R}^n$, input $\bu(t) \in \mathbb{R}^k$, linear operator $\mathbf{A} \in \mathbb{R}^{n \times n}$, quadratic operator $\mathbf{H} \in \mathbb{R}^{n \times n^2}$, input operator $\mathbf{B} \in \mathbb{R}^{n \times k}$ and constant term $\mathbf{c} \in \mathbb{R}^{n \times 1}$. The finite-dimensional system in \eqref{quad} comprises a prototypical structure that is typically obtained after the spatial discretization of a PDE with a quadratic nonlinear term. 

We note that the term $\bx(t) \otimes \bx(t)$ comprises of $n^2$ repeating terms, (e.g. $\bx_i(t) \bx_j(t)$ and $\bx_j(t) \bx_i(t)$, where $i \neq j$), only ${n \choose 2}$ of which are unique. 

Matrix $\bH$ can be written with respect to $n$ sub-matrices, as

\begin{equation}
\label{sub-matrices}
\bH = \begin{bmatrix} {\bH}_1, \ldots, {\bH}_n \end{bmatrix}
\end{equation}

\noindent such that the element of matrix $\bH$ in the $j$-th row and $k$-th column of the $i$-th sub-matrix is denoted by three indices, $h_{i_{j,k}}$. This notation will be used in the following.

We will consider exclusively systems with an energy-preserving $\bH$. Such systems are examined in \cite{Schlegel2015} and are shown to satisfy

\begin{equation}
\label{energy}
h_{i_{j,k}}+h_{i_{k,j}}+h_{j_{i,k}}+h_{j_{k,i}}+h_{k_{i,j}}+h_{k_{j,i}}=0.
\end{equation}

\noindent If the quadratic term is energy-preserving, it does not contribute to the overall kinetic energy of the system, but redistributes the energy among the system DOFs. 
 In \cite{goyal2023guaranteed}, it has been shown that \eqref{energy} can equivalently be written as 
 \begin{equation}
 	\label{energy_pres}
 	\bx(t)^\top\mathbf{H} \left( \bx(t) \otimes \bx(t) \right) =0, \; \forall \bx(t) \in \mathbb{R}^n.
 \end{equation}
 Equation \eqref{energy} introduces $n(n^2+5)/6$ constraints to the entries of matrix $\mathbf{H}$. The question we aim to answer is whether the linear constraints in \eqref{energy} (and thus, equivalently \eqref{energy_pres}) can be a priori encoded to $\bH$, via an equivalent representation, denoted as $\tilde{\mathbf{H}}$. The equivalence of operators $\bH$ and $\tilde{\mathbf{H}}$ is defined as

\begin{equation}
\label{equiv}
\mathbf{H} \left( \bx(t) \otimes \bx(t) \right) =\tilde{\mathbf{H}} \left( \bx(t) \otimes \bx(t) \right), \; \forall \bx(t) \in \mathbb{R}^n,
\end{equation}

\noindent such that $\bH$ can be substituted by $\tilde{\mathbf{H}}$ in \eqref{quad}. It is useful to employ the previously introduced notation and express \eqref{equiv} in terms of quadratic entry summands of permuted $i,j$ indices as

\begin{equation}
\label{dualperm}
\tilde{h}_{i_{k,j}}+\tilde{h}_{j_{k,i}}=h_{i_{k,j}}+h_{j_{k,i}}.
\end{equation}

\noindent In \cite{goyal2023guaranteed}, a parameterization $\tilde{\bH}$ such that

\begin{equation}
\label{skew}
\tilde{\bH} = \begin{bmatrix} \tilde{\bH}_1, \ldots, \tilde{\bH}_n \end{bmatrix}, \qquad \tilde{\mathbf{H}}_i= -{\tilde{\mathbf{H}}_i}^\top, \; i=\left[1, \dots, n \right],
\end{equation} 

\noindent is proposed, such that $\tilde{\bH}$ is energy-preserving. By writing \eqref{skew} using the introduced notation, as

\begin{equation}
\label{skew1}
\tilde{h}_{i_{k,j}}=-\tilde{h}_{i_{j,k}},
\end{equation}

\noindent the aforementioned result can be understood by substituting \eqref{skew1} to \eqref{energy} for $\tilde{\mathbf{H}}$. 

In this work, we will make use of representation \eqref{skew}. In particular, we answer whether \textit{every} energy-preserving nonlinearity $\mathbf{H}$ admits the representation in \eqref{skew}. Making use of the quadratic operators equivalence in \eqref{equiv}, we reformulate the question as to whether there exists an equivalent matrix $\tilde{\mathbf{H}}$ for \textit{any} energy-preserving $\mathbf{H}$, such that $\tilde{\mathbf{H}}$ satisfies \eqref{skew}. Based on that question, we formulate the main result of this work.

\begin{theorem*} For any energy-preserving quadratic matrix $\mathbf{H} \in \mathbb{R}^{n \times n^2}$, such that $\mathbf{x}^\top\mathbf{H} \left( \mathbf{x} \otimes \mathbf{x} \right)=0, \; \forall \mathbf{x} \in \mathbb{R}^n$, there always exists an equivalent $\tilde{\mathbf{H}}$ with a structure
\begin{equation}
\tilde{\bH} = \begin{bmatrix} \tilde{\bH}_1, \ldots, \tilde{\bH}_n \end{bmatrix}, \qquad \tilde{\mathbf{H}}_i=-\tilde{\mathbf{H}}^\top_i, \; \; i=\left[1, \dots, n \right],
\end{equation}
such that $\mathbf{H} \left( \mathbf{x} \otimes \mathbf{x} \right)=\tilde{\mathbf{H}} \left( \mathbf{x} \otimes \mathbf{x} \right)$. Some $n^3/6-n^2/2+n/3$ entries of $\tilde{\mathbf{H}}$ can be freely selected. 
\end{theorem*}
\label{theor1}

\begin{proof}
To prove this, we start by enumerating all linearly independent equations necessary to define $\tilde{\mathbf{H}}$, given $\mathbf{H}$. We first focus on the equivalence requirement in \eqref{dualperm}. We classify the index configurations in three cases: $(a): \; \{i,j,k \; \rvert \; i=j \neq k \; \cup \; i=k \neq j\; \cup \; j=k \neq i\}$, $(b): \;  \{i,j,k \; \rvert \; i \neq j \neq k\}$, $(c): \;  \{i,j,k \; \rvert \; i=j=k\}$.

\noindent $(a)$: We consider $k>i=j$ for equation \eqref{dualperm}. Then we get the following $n(n-1)/2$ equations:

\begin{equation}
\label{dualp1}
\tilde{h}_{i_{k,i}}=h_{i_{k,i}}.
\end{equation}

\noindent Conversely, we consider $j=k>i$. By skew-symmetry in $\tilde{\mathbf{H}}_i$ and energy-preservation for $j=k$, we get $n(n-1)/2$ equations 

\begin{equation}
\label{dualp2}
\tilde{h}_{k_{k,i}}=-h_{k_{i,k}}-\cancelto{0}{\tilde{h}_{i_{k,k}}}.
\end{equation}

\noindent Up to now, we have obtained $n(n-1)$, linearly independent equations,  uniquely defining elements of $\tilde{\mathbf{H}}$. Having set $k>i$, these equations refer to elements of $\tilde{\mathbf{H}}$ below the diagonal of each sub-matrix, with one repeating index. Thus, they are linearly independent from the equations for skew-symmetry in \eqref{skew1}.

\noindent Resulting equations for $i=k\neq j$ are

\begin{equation}
\tilde{h}_{i_{i,j}}=-h_{i_{j,i}}.
\end{equation}

\noindent However, these equations are all redundant. For $i>j$ they coincide with \eqref{dualp2}. For $i<j$, they are already satisfied by \eqref{dualp1} and skew-symmetry \eqref{skew1}.

\noindent $(b)$: We continue by considering $i \neq j \neq k$. This results to $n(n-1)(n-2)/2$ unique equations from \eqref{dualperm}. We notice that for a triplet $\{(i,j,k) \rvert i \neq j \neq k \}$ one needs to impose only two out of the three equations \eqref{dualperm}. Due to energy preservation for matrix $\mathbf{H}$ and skew-symmetry of $\tilde{\mathbf{H}}$, we get:

\begin{multline}
\label{rel}
\begin{aligned}
h_{i_{j,k}}+h_{i_{k,j}}+h_{j_{i,k}}+h_{j_{k,i}}+h_{k_{i,j}}+h_{k_{j,i}} &= 
\tilde{h}_{i_{j,k}}+\tilde{h}_{i_{k,j}}+\tilde{h}_{j_{i,k}}+\tilde{h}_{j_{k,i}}+\tilde{h}_{k_{i,j}}+\tilde{h}_{k_{j,i}} &=0.
\end{aligned}
\end{multline}

\noindent Thus, if for every triplet $\{(i,j,k) \rvert i \neq j \neq k \}$, 
two out of the three equations in \eqref{dualperm} are satisfied, the third one will also be satisfied due to \eqref{rel}.

\noindent We thus enumerate $n(n-1)(n-2)/3$ linearly independent equations for $i \neq j \neq k$:

\begin{equation}
\label{dualp4}
\tilde{h}_{i_{k,j}}=h_{i_{k,j}}+h_{j_{k,i}}-\tilde{h}_{j_{k,i}}.
\end{equation}

\noindent These equations are linearly independent from the equations in \eqref{dualp1}, \eqref{dualp2} and their corresponding skew-symmetry conditions, since index sets $(a)$ and $(b)$ do not intersect. What remains is to show that these equations are linearly independent from the corresponding skew-symmetry equations. \eqref{dualp4} correlates entries in the same row of $\tilde{\mathbf{H}}$, while skew-symmetry equations correlate entries in the same sub-matrix of $\tilde{\mathbf{H}}$. Thus, using the two equations for $i \neq j \neq k$ cannot result to linear dependencies (since $\tilde{h}_{i_{k,j}}$ cannot be expressed with respect to $\tilde{h}_{i_{j,k}}$). Thus all $n(n-1)(n-2)/3$ equations in \eqref{dualp4} are also linearly independent from all other previously enumerated equations.

\noindent $(c)$: The case for $i=j=k$ is redundant, due to skew-symmetry conditions. 

From $(a)$, $(b)$ and $(c)$, we count $n(n-1)+n(n-1)(n-2)/3$ linearly independent equations originating from all dual permutations in \eqref{dualperm}, i.e. $n(n-1)+n(n-1)(n-2)/3$ linear constraints on $\tilde{\mathbf{H}}$, given $\mathbf{H}$. The skew-symmetry of the sub-matrices in $\tilde{\mathbf{H}}$ (eq. \eqref{skew1}) introduces another $n^2(n+1)/2$ equations, which, as we showed, are also linearly independent from those in \eqref{dualp1}, \eqref{dualp2}, \eqref{dualp4}. 

In total, we have obtained a system of $n(n-1)+n(n-1)(n-2)/3+n^2(n+1)/2=5/6n^3+n^2/2-n/3$ linearly independent equations for the $n^3$ unknowns for the new matrix, $\tilde{\mathbf{H}}$. We compute that $5/6n^3+n^2/2-n/3=n^3$ for $n=1,2$ and $5/6n^3+n^2/2-n/3<n^3$ for $n>3$. This means that if the (possibly unknown) quadratic nonlinear terms are energy preserving, there always exists a representation $\tilde{\mathbf{H}}$ of the quadratic term with antisymmetric sub-matrices, such that $\mathbf{H} \left( \mathbf{x} \otimes \mathbf{x} \right)=\tilde{\mathbf{H}} \left( \mathbf{x} \otimes \mathbf{x} \right), \; \forall \mathbf{x} \in \mathbb{R}^n$.  $n^3/6-n^2/2+n/3$ entries of $\tilde{\mathbf{H}}$ for $i\neq j\neq k$ can be chosen freely.
\end{proof}

Using matrix notation, an alternative presentation of this proof is given in Appendix \ref{appendix}. The above proof can be straightforwardly recast for other representations $\tilde{\mathbf{H}}$, other than that of skew-symmetric sub-matrices. For example, it might be interesting to use a representation that encodes skew-symmetry between rows of the different sub-matrices of $\tilde{\mathbf{H}}$, written as

\begin{equation}
\label{fact2}
\tilde{h}_{i_{k,j}}=-\tilde{h}_{k_{i,j}}.
\end{equation}

\noindent  A key of the above argumentation for representations such as skew-symmetry in \eqref{skew1} or the row-wise ``skew-symmetry" in \eqref{fact2} is that, since \eqref{rel} holds, only $2/3$ of the equations for $i \neq j \neq k$ need to be enforced for the corresponding dual permutations in \eqref{dualperm}. More complex representations than \eqref{fact2} could also be constructed, as long as there is a nullspace to the system of linearly independent equations defining $\tilde{\mathbf{H}}$. We will not pursue this task further in this work.

\section{Algorithm for quadratic matrix representation}
\label{sec:algo}

Based on the presented proof, we construct \Cref{alg:cap}, which transforms any energy-preserving $\mathbf{H}$ to an equivalent $\tilde{\mathbf{H}}$ with skew-symmetric sub-matrices. This procedure also showcases the row echelon form of the system of $5/6n^3+n^2/2-n/3$ linearly independent equations.

\begin{algorithm}
\caption{Representation to energy-preserving matrix with skew-symmetric sub-matrices.}\label{alg:cap}
\begin{algorithmic}
 \State \textbf{Input}: Energy-preserving quadratic matrix $\mathbf{H}$
 \State \textbf{Output}: Equivalent quadratic matrix $\tilde{\mathbf{H}}$, with skew-symmetric sub-matrices.
\begin{enumerate}
 \item Set all sub-matrix diagonal entries of $\tilde{\mathbf{H}}$ to zero (due to \eqref{skew1}). 
 \item Compute $\tilde{\mathbf{H}}$ entries in \eqref{dualp1} and \eqref{dualp2}. 
 \item Apply corresponding skew-symmetry constraints \eqref{skew1}.
 \item Set $n^3/6-n^2/2+n/3$ entries of $\tilde{\mathbf{H}}$ for $i \neq j \neq k$, at most one per triplet $\{i,j,k\}$. 
 \item For all entries in 4, apply equations \eqref{dualp4} and \eqref{skew1} iteratively, until matrix completion.
 \end{enumerate}
 \label{algorithm1}
\end{algorithmic}
\end{algorithm}

We illustrate the above algorithm for an arbitrary, energy-preserving matrix $\mathbf{H} \in \mathbb{R}^{3 \times 9}$:

\begin{figure}[H]
\begin{subfigure}{.45\textwidth}
  \centering
  \includegraphics[width=.80\linewidth]{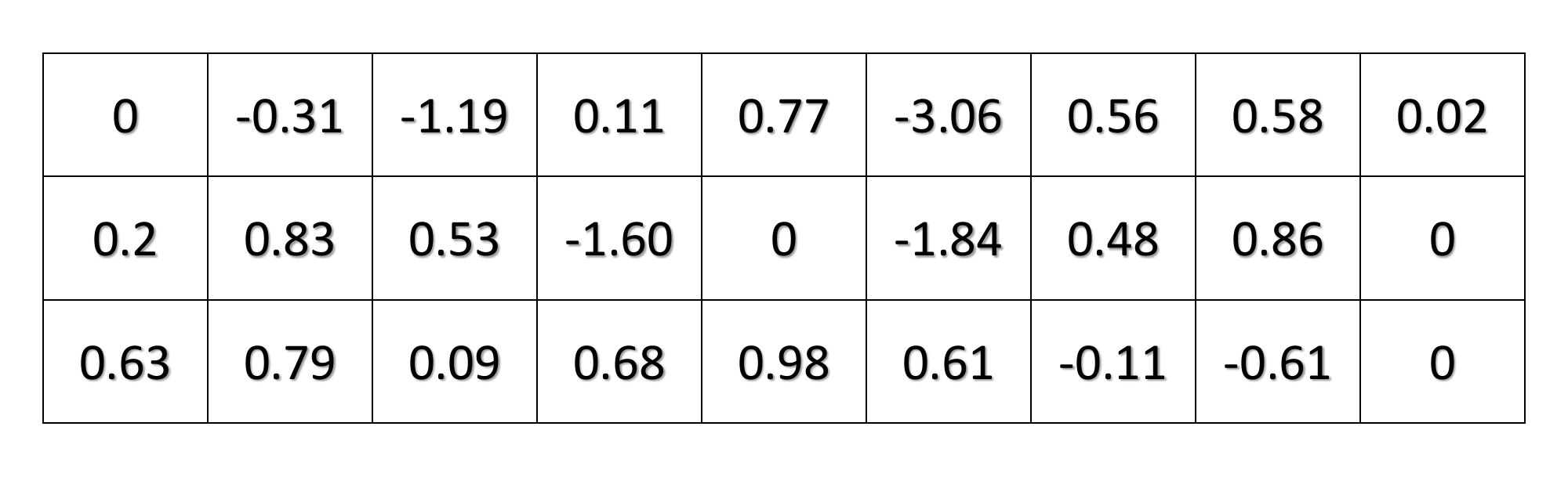}
  \caption{Arbitrary, energy-preserving $\mathbf{H} \in \mathbb{R}^{3 \times 9}$. It can be checked that indeed $\mathbf{x}^\top\mathbf{H} \left( \mathbf{x} \otimes \mathbf{x} \right)=0, \; \forall \mathbf{x} \in \mathbb{R}^3$.}
  \label{fig:sub1}
\end{subfigure}%
\hfill
\begin{subfigure}{.45\textwidth}
  \centering
  \includegraphics[width=.80\linewidth]{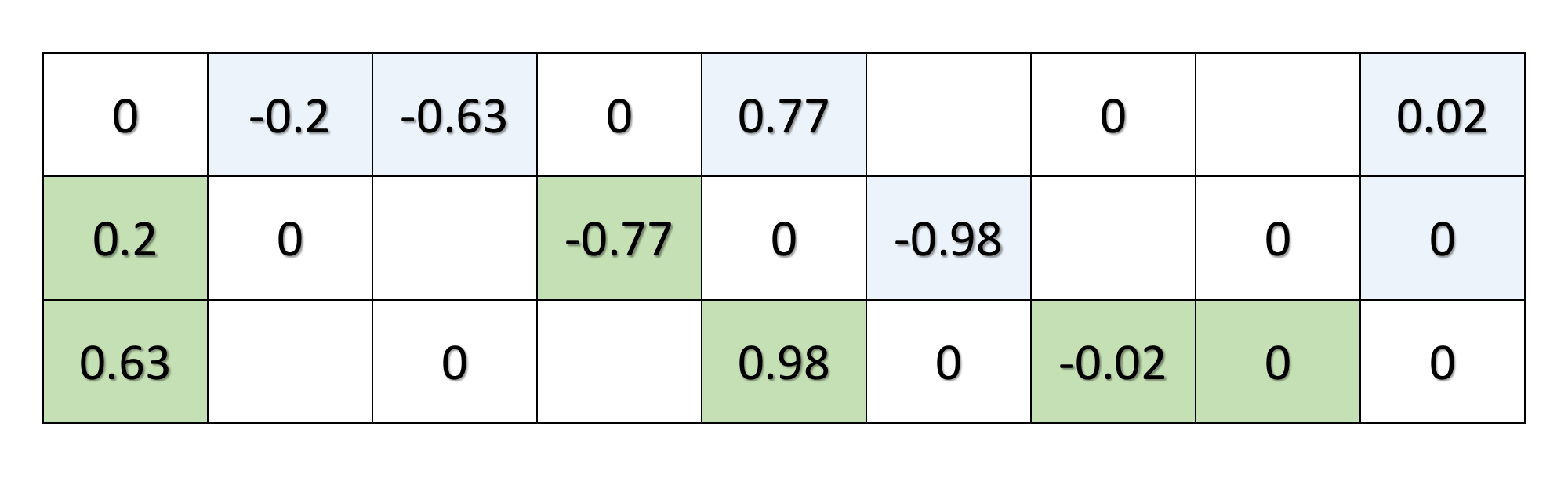}
  \caption{ Apply steps 1 to 3 of \Cref{alg:cap}: Entries from \eqref{dualp1} and \eqref{dualp2} are marked in green, while entries from \eqref{skew1} are marked in blue.}
  \label{fig:sub2}
  \end{subfigure}
\label{fig:test}

\bigskip 
\begin{subfigure}{.45\textwidth}
  \centering
  \includegraphics[width=.80\linewidth]{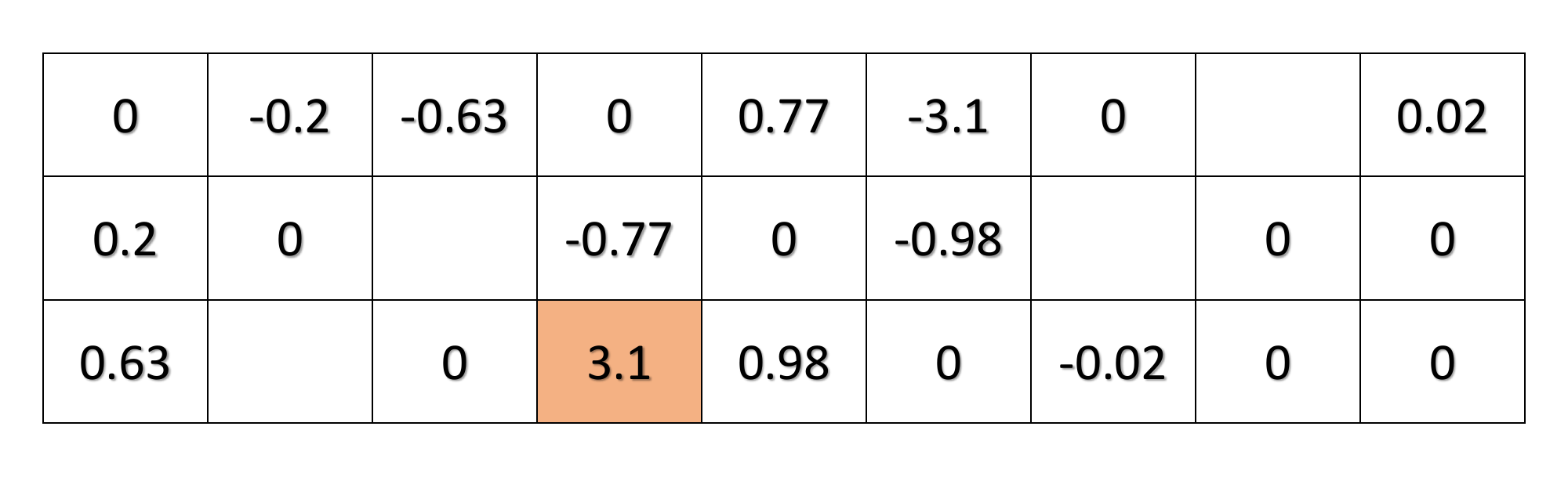}
  \caption{For $n=3$, we get $n^3/6-n^2/2+n/3=1$ freely selected entries (marked in red). We set entry $\tilde{h}_{2_{3,1}}=3.1$ (step 4 of \Cref{alg:cap}).}
  \label{fig:sub3}
\end{subfigure}%
\hfill
\begin{subfigure}{.45\textwidth}
  \centering
  \includegraphics[width=.80\linewidth]{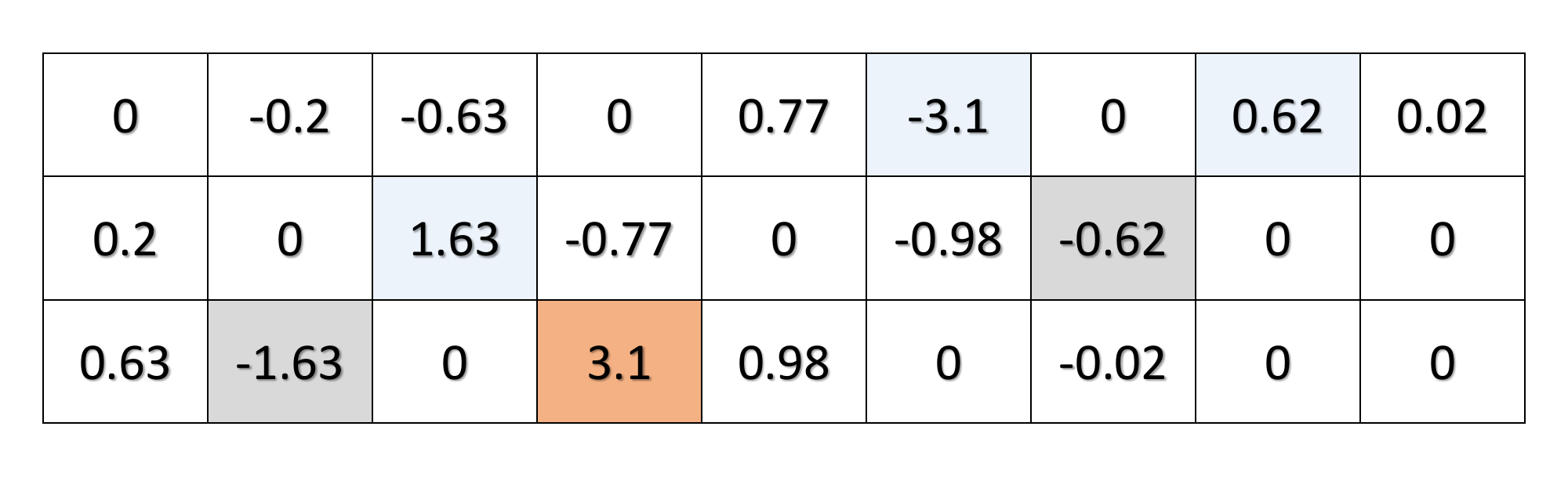}
  \caption{Apply step 5 of \Cref{alg:cap} for entries in \eqref{dualp4}, marked in grey. Skew-symmetric entries are marked in blue.}
  \end{subfigure}
  \label{fig:subfinal}
\caption{Illustration of \Cref{alg:cap} for an arbitrary, energy-preserving matrix $\mathbf{H} \in \mathbb{R}^{3 \times 9}$.}
\end{figure}

\noindent We observe that the resulting quadratic matrix $\tilde{\mathbf{H}}$ has skew-symmetric blocks and is thus energy-preserving. One can check that indeed $\tilde{\mathbf{H}} \left( \mathbf{x} \otimes \mathbf{x} \right)= \mathbf{H} \left( \mathbf{x} \otimes \mathbf{x} \right)$. 

\begin{figure}[!htb]
  \centering
  \includegraphics[width=.38\linewidth]{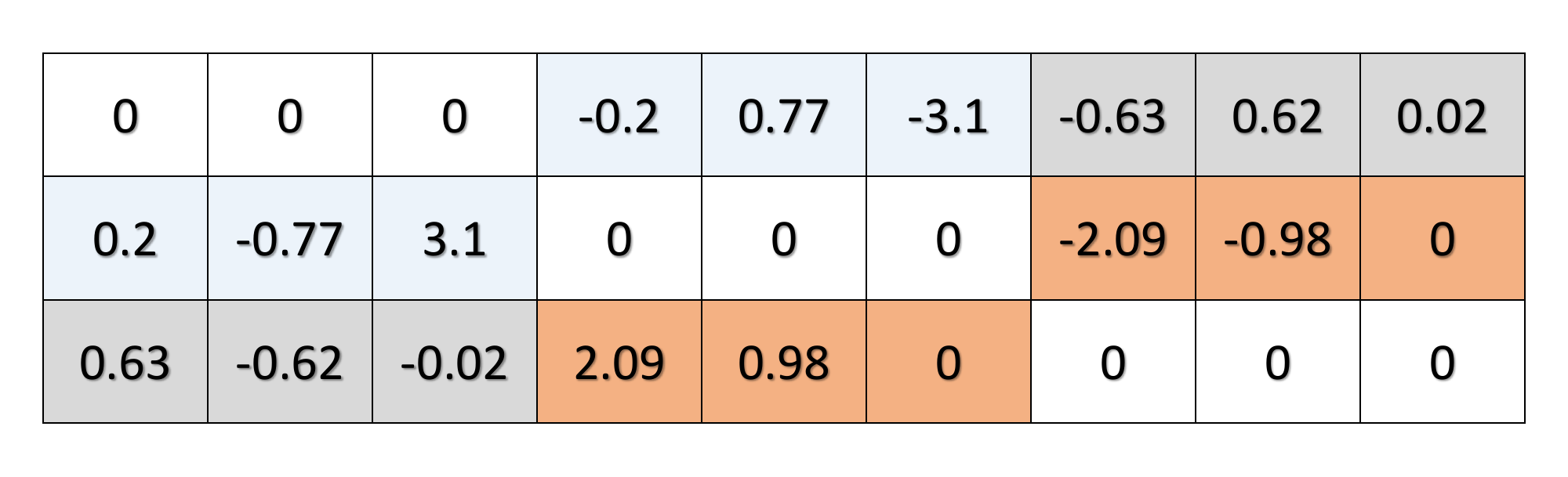}
  \caption{Alternative representation via \eqref{fact2}. Note that both this representation as well as the one in \Cref{fig:subfinal} leave the entries $h_{i_{j,i}}, \forall i,j \in \left[1, \dots, n\right]$ in \Cref{fig:sub1} unchanged. These correspond to non-repeating terms in $\mathbf{H}\left( \mathbf{x} \otimes \mathbf{x} \right)$.}
  \label{fig:alt}
\end{figure}

\noindent Following representation \eqref{fact2}, the operator in \Cref{fig:sub1} could be equivalently represented as shown in \Cref{fig:alt}. The ``skew-symmetric" feature of this representation concerns the rows of each sub-matrix of  $\tilde{\mathbf{H}}$. However, we will focus on representation \eqref{skew1} for the remainder of this work.

\section{Application to Operator Inference}
\label{sec:OpInf}

\subsection{The OpInf framework}
In this section we use the proven theorem in the context of physics-informed, non-intrusive modeling and in particular OpInf \cite{Peherstorfer2016} and propose a new approach for guaranteeing energy preservation of the quadratic term via a sequence of closed-form LS problem solutions.

We first briefly present the OpInf framework for physics-informed, non-intrusive model reduction. We focus on quadratic systems, which emerge from many physics and engineering applications \cite{Benner2021, Qian2020, Kalur2021, McQuarrie2021}. For a given discretized system with state vector $\mathbf{x}(t) \in \mathbb{R}^n$ and input vector $\mathbf{u}(t) \in \mathbb{R}^k$, assuming state access at timesteps $[ t_1, ..., t_m]$, we construct a snapshot matrix $\mathbf{X} \in \mathbb{R}^{n \times m}$ with state data and a snapshot matrix $\mathbf{U}$ with input data as

\begin{equation}
    \mathbf{X}=\begin{bmatrix}
    \mathbf{x}(t_1) & ...&   \mathbf{x}(t_m) 
    \end{bmatrix}, \qquad
    \mathbf{U}=\begin{bmatrix}
    \mathbf{u}(t_1) & ...&   \mathbf{u}(t_m) 
    \end{bmatrix}. 
\end{equation}

Performing the Singular Value Decomposition (SVD) of $\mathbf{X}$ and projecting $\mathbf{X}$ to its leading $r$ SVD modes $V \in \mathbb{R}^{n \times r}$, leads to a reduced-order snapshot matrix $\hat{\mathbf{X}}$. We also assume time derivative data is available or can be numerically computed and assembled in a corresponding snapshot matrix $\dot{\hat{\mathbf{X}}}$. Within the framework of OpInf, we aim to infer a ROM of quadratic model from reduced state data $\hat{\mathbf{X}}$, input data $\mathbf{U}$ and corresponding time derivative data $\dot{\hat{\mathbf{X}}}$, by solving a linear LS problem:

\begin{equation}
\label{OpInf}
\min_{\hat{\mathbf{c}},\hat{\mathbf{A}},\hat{\mathbf{H}},\hat{\mathbf{B}}}{ \left\| \hat{\mathbf{c}}\mathbf{1}^\top_{r\times m} + \hat{\mathbf{A}} \hat{\mathbf{X}} + \hat{\mathbf{H}} ( \hat{\mathbf{X}} \otimes \hat{\mathbf{X}}) + \hat{\mathbf{B}} \mathbf{U} - \dot{\hat{\mathbf{X}}} \right\|}_F^2.
\end{equation}

\noindent One of the advantages of OpInf is that \eqref{OpInf} admits a closed-form solution. \eqref{OpInf} can be recast to a series of LS problems for every  $i \in [1, ..., r]$, as presented in the original contribution of the method \cite{Peherstorfer2016}. Specifically, the minimization problem in \eqref{OpInf} is equivalent to

\begin{equation}
\label{seqls}
\sum_{j=1}^r{\min_{\hat{c}_j,\hat{\mathbf{A}}_{j,:},\hat{\mathbf{H}}_{j,:},\hat{\mathbf{B}}_{j,:}}{ \left\| \hat{c}_j\mathbf{1}_{1\times m} + \hat{\mathbf{A}}_{j,:} \hat{\mathbf{X}} + \hat{\mathbf{H}}_{j,:} ( \hat{\mathbf{X}} \otimes \hat{\mathbf{X}}) + \hat{\mathbf{B}}_{j,:} \mathbf{U} - \dot{\hat{\mathbf{X}}}_{j,:} \right\|}_2^2}.
\end{equation}

\noindent Thus \eqref{OpInf} can be split to $r$ LS problems

\begin{equation}
\label{seqls2}
\min_{\hat{c}_j,\hat{\mathbf{A}}_{j,:},\hat{\mathbf{H}}_{j,:},\hat{\mathbf{B}}_{j,:}}{ \left\| \hat{c}_j\mathbf{1}_{1 \times m} + \hat{\mathbf{A}}_{j,:} \hat{\mathbf{X}} + \hat{\mathbf{H}}_{j,:} ( \hat{\mathbf{X}} \otimes \hat{\mathbf{X}}) + \hat{\mathbf{B}}_{j,:} \mathbf{U} - \dot{\hat{\mathbf{X}}}_{j,:} \right\|}_2^2, \; j=[1, \dots, r].
\end{equation}

\noindent Concatenating the different terms and retaining only the $r(r+1)/2$ unique quadratic term entries, \eqref{seqls2} can be rewritten following the notation in \cite{Peherstorfer2016}, as

\begin{equation}
\label{seqls3}
\min_{\bm{\mathcal{O}}_{j,:}}{ \left\| \bm{\mathcal{O}}_{j,:} \bm{\mathcal{D}} - \dot{\hat{\mathbf{X}}}_{j,:} \right\|}_2^2, \; j=[1, \dots, r],
\end{equation}

\noindent with $\bm{\mathcal{O}}=[\hat{\mathbf{c}} \; \hat{\mathbf{A}} \; \hat{\mathbf{H}} \; \hat{\mathbf{B}}] \in \mathbb{R}^{r \times (1+r+r(r+1)/2+k)}$, $\bm{\mathcal{D}}=[\mathbf{1}_{1 \times m}^\top \; \hat{\mathbf{X}}^\top \; {\left(\hat{\mathbf{X}}^2\right)}^\top \; \hat{\mathbf{U}}^\top]^\top \in \mathbb{R}^{(1+r+r(r+1)/2+k) \times m}$.
Each LS problem in \eqref{seqls3} has a closed-form solution, which is related to the dynamics of an SVD mode $j$, where an increasing $j$ indicates a monotonically decreasing significance in terms of system energy. In practice, \eqref{seqls3} should be complemented with a -typically $L_2$- regularization term. This term penalizes the LS solution norm with a factor $\lambda$, leveraging the LS error and the norm of the obtained solution. The optimal $\lambda$ value can be selected via an L-curve criterion, see \cite{McQuarrie2021, Kramer2019}. We note that formulations \eqref{OpInf} and \eqref{seqls} are equivalent only if a common regularization value $\lambda$ is selected for all $r$ LS problems in \eqref{seqls3} as in \eqref{OpInf}. Then, the $L_2$ norm minimization for all summands in \eqref{seqls} is equivalent to the Frobenius norm minimization in \eqref{OpInf}. However, formulation \eqref{seqls3} allows for more flexibility in choosing the optimal $\lambda$ value independently for each of the $r$ LS problems. In the following, we will denote formulation \eqref{seqls3} as \texttt{Seq\_OpInf}.

\subsection{Sequential Operator Inference formulation for energy-preserving nonlinearities}

The successful application of \eqref{seqls3} is significant for physics-informed data-driven inference, where only experimental or proprietary software simulation data are available for the system states.
In the direction of physics-informed constraints, the proven representation for quadratic operators in \Cref{sec:fact} could be employed in applications where the quadratic operator $\mathbf{H}$ of the system at hand is known to be (exactly or approximately) energy preserving, i.e. $\mathbf{x}^\top\mathbf{H} \left( \mathbf{x} \otimes \mathbf{x} \right) \approx 0$. Having proven the admissibility of representation \eqref{skew} for any energy-preserving quadratic nonlinearity, we employ formulation \eqref{seqls2} to suggest a closed-form solution for the inference of systems with energy-preserving quadratic nonlinearities.

The sequential solution of \eqref{seqls2} was successfully employed to enforce the intrusive operators property of nested-ness to the inferred operators \cite{NAretz}. In a similar fashion, we enforce the skew-symmetric sub-matrix representation \eqref{skew}, proven in \Cref{sec:fact}, for cases of energy-preserving quadratic operators to $\hat{\mathbf{H}}$.

We use \eqref{seqls2} since we need all $r^2$ terms of the quadratic operator to enforce \eqref{skew1}, instead of the unique $r(r+1)/2$ terms used in \texttt{Seq\_OpInf}.

We will use notation $h_{i_{j,k}}$ to illustrate the process of enforcing \eqref{skew1} to \eqref{seqls2}. First, we solve \eqref{seqls2} e.g. for $j=1$, neglecting the contribution of $h_{i_{1,1}}$, (since it is identically zero due to skew-symmetry). The $r^2-r$ inferred terms of the first row of $\hat{\mathbf{H}}$ dictate the values of all entries $h_{i_{j,k}}$ with $k=1$, i.e. the first column of each sub-matrix of $\hat{\mathbf{H}}$, due to skew-symmetry \eqref{skew1}. Then, for $j=2$, $h_{i_{2,2}}$ terms are neglected. Hence, for $j=2$, $2r$ entries of $\hat{\mathbf{H}}_{2,:}$ are fixed due to sub-matrix skew-symmetry conditions: $r$ sub-matrix diagonal entries are identically set to $0$ and $r$ entries ($k=1$) are fixed due to the LS solution for $j=1$. The quadratic terms corresponding to those $r$ entries for $k=1$ can be moved to the right-hand side of the $2$nd LS problem. Hence, the remaining $r^2-2r$ entries for $j=2$ will be inferred by the $2$nd LS problem. The process continues by solving problems \eqref{seqls2} through $j=r$, with the $j$-th LS problem involving the inference of $r^2-jr$ quadratic term entries: $r$ sub-matrix diagonal entries are a priori excluded from inference (identically zero) and $(j-1)r$ entries are fixed due to previous LS solutions.

We denote as $f_j$, the subset of $(j-1)r$ fixed quadratic operator entries for the $j$-th LS problem, due to sub-matrix skew-symmetry. The remaining $r^2-jr$ entries of $\hat{\mathbf{H}}_{j,:}$ (sub-matrix diagonals excluded) that need to be inferred are denoted as $d_j$. Then the above approach can be written formally as 

\begin{multline}
\label{seqls4}
\min_{\hat{c}_j,\hat{\mathbf{A}}_{j,:},\hat{\mathbf{H}}_{j,d_j},\hat{\mathbf{B}}_{j,:}}{ \left\| \hat{c}_j\mathbf{1}^\top_r + \hat{\mathbf{A}}_{j,:} \hat{\mathbf{X}} + \hat{\mathbf{H}}_{j,d_j} ( \hat{\mathbf{X}} \otimes \hat{\mathbf{X}})_{d_j,:} +\hat{\mathbf{B}}_{j,:} \mathbf{U} - (\dot{\hat{\mathbf{X}}} - \hat{\mathbf{H}}_{j,f_j} ( \hat{\mathbf{X}} \otimes \hat{\mathbf{X}}))_{f_j,:} \right\|}_F^2, \; j=[1, \dots, r].
\end{multline}

Similarly to \texttt{Seq\_OpInf} in \eqref{seqls3}, we rewrite \eqref{seqls4} in more compact notation as

\begin{equation}
\label{seqls5}
\min_{\bm{\mathcal{O}}_{j,:}}{ \left\| \bm{\mathcal{O}}_{j,:} \bm{\mathcal{D}}_{j,:} - {\bm{\mathcal{F}}}_{j,:} \right\|}_2^2, \; j=[1, \dots, r],
\end{equation}

\noindent with $\bm{\mathcal{O}}_{j,:}=[\hat{\mathbf{c}}_j \; \hat{\mathbf{A}}_{j,:} \; \hat{\mathbf{H}}_{j,d_j} \; \hat{\mathbf{B}}_{j,:}] \in \mathbb{R}^{1 \times (1+r+r^2-jr+k)}$, $\bm{\mathcal{D}}_{j,:}=[\mathbf{1}_{1 \times m}^\top \; \hat{\mathbf{X}}^\top \; {( \hat{\mathbf{X}} \otimes \hat{\mathbf{X}})_{d_j,:}}^\top \; \hat{\mathbf{U}}^\top]^\top \in \mathbb{R}^{(1+r+r^2-jr+k) \times m}$ and $\bm{\mathcal{F}}_{j,:} = \dot{\hat{\mathbf{X}}}_{j,:} - \hat{\mathbf{H}}_{j,f_j} ( \hat{\mathbf{X}} \otimes \hat{\mathbf{X}}))_{f_j,:} \in \mathbb{R}^{1 \times m}$.

We will denote this formulation as \texttt{Seq\_OpInf\_EP}. The $j$-th LS problem in \eqref{seqls5} has a closed-form solution, sequentially depending on all preceding LS solutions. Given that the involved LS problems in \texttt{Seq\_OpInf} and \texttt{Seq\_OpInf\_EP} are solved in the reduced dimension $r$, the additional clock time of sequentially solving the LS problems is negligible. Observing \eqref{seqls3} and \eqref{seqls5}, it should also be noted that configuring \texttt{Seq\_OpInf} to \texttt{Seq\_OpInf\_EP} requires negligible programming effort. In \texttt{Seq\_OpInf\_EP}, the total quadratic matrix entries to be inferred are $r^2(r-1)/2$ (less than $r^2(r+1)/2$ in the general case of \texttt{Seq\_OpInf}) with each LS problem involving $r$ less terms than the previous one. Both \texttt{Seq\_OpInf} in \eqref{seqls3} and \texttt{Seq\_OpInf\_EP} in \eqref{seqls5} are, in practice, complemented by an $L_2$ regularization term. Formulation \texttt{Seq\_OpInf\_EP} in \eqref{seqls5} entails a sequential dependence in LS solutions and corresponding regularization values: The quality of each LS solution will be affected by that of all preceding solutions. 

\subsection{Numerical Example: 2D Burgers' equation}

\noindent The feasibility of \texttt{Seq\_OpInf\_EP}, given the proven theorem in \Cref{sec:fact}, is showcased for a simple 2D benchmark. We simulate the viscous, 2D Burgers' equation

\begin{equation}
    \label{burg}
\frac{\partial u}{\partial t}=c u  \nabla \cdot u + \nu \nabla^2{u},
\end{equation}

\noindent for $(x,y) \in [0, 1]^2$, $c=0.2$, $\nu=2 \times 10^{-3}$, periodic boundary conditions and initial condition $u(t=0)=cos(2\pi x)cos(2\pi y)$. For this problem, we expect that the quadratic nonlinear term $u  \nabla \cdot u$ is indeed energy preserving, since

\begin{equation}
\label{integral}
\int_{0}^{1} \int_{0}^{1} u^2  \nabla \cdot u \; dx \; dy =0,
\end{equation}

\noindent due to periodic boundary conditions. 

\noindent We simulate this system over $T=4 s$, with a timestep $\Delta t =0.01 \; s$, a second-order discretization in space and $\Delta x =0.02 \; m$. After discretizing in space, \eqref{burg} becomes

\begin{equation}
\label{burg2}
\dot{\bx}(t) =\mathbf{A}_v \bx(t) + \mathbf{H}_a \left( \bx(t) \otimes \bx(t) \right),
\end{equation}

\noindent where $\mathbf{A}_v \bx(t)$ corresponds to the discretized viscous, linear dissipation term and $\mathbf{H}_a \left( \bx(t) \otimes \bx(t) \right)$ corresponds to the discretized advection term. Given an appropriate discretization scheme, $\mathbf{H}_a$ introduces a quadratic, energy-preserving nonlinearity, since \eqref{integral} is the continuous counterpart of \eqref{energy_pres} for $\mathbf{H}_a$.

In the OpInf context, we aim to infer a quadratic ROM (i.e. \eqref{quad} with $\mathbf{x} \in \mathbb{R}^r$, $\mathbf{B}=\mathbf{0}$ and $\mathbf{c}=\mathbf{0}$) from simulation data. Essentially we aim to infer operators $\mathbf{A} \in \mathbb{R}^{r \times r}$ and $\mathbf{H} \in \mathbb{R}^{r \times r^2}$. For this task, both standard \texttt{Seq\_OpInf} \eqref{seqls3} and \texttt{Seq\_OpInf\_EP} in \eqref{seqls5} (with the physics-based knowledge of $\mathbf{x}^\top\mathbf{H} \left( \mathbf{x} \otimes \mathbf{x} \right) = 0$) are employed.

\noindent Since we aim to only compare the two approaches, we use all available data as training data. For both \texttt{Seq\_OpInf} and \texttt{Seq\_OpInf\_EP}, we complement the LS problem (\eqref{seqls3} and \eqref{seqls5} accordingly) with an $L_2$ regularization term, scaled accordingly for the linear and quadratic operator entries. This term penalizes the $L_2$ norm of the LS solution with a factor $\lambda$, for which $50$ logarithmically spaced values in the interval $[ 10^{-5}, 10^{3}]$ are tested. The quadratic term entries are regularized with factor $r \lambda$, to account for the different scaling of the linear and quadratic terms under different ROM dimensions $r$. The optimal $\lambda$ regularization factor is selected via an $L$-curve criterion in both cases.

\begin{figure}[H]
\begin{subfigure}[t]{.50\textwidth}
  \centering
  \includegraphics[width=1.0\linewidth]{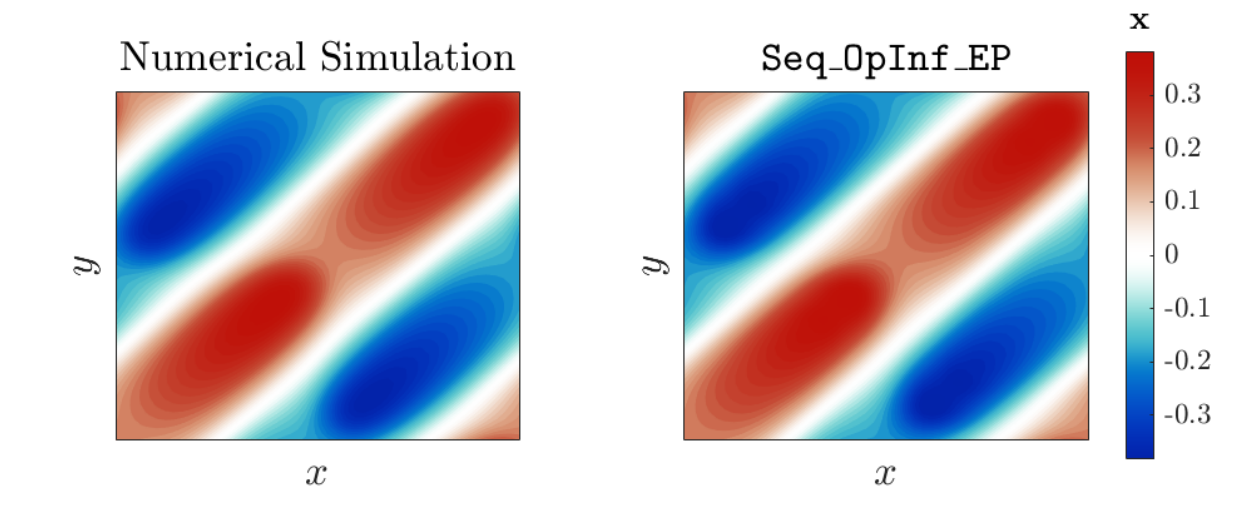}
  \caption{State prediction for 2D Burgers' in \eqref{burg}, at $t=4 \;s$ and $r=15$, with \texttt{Seq\_OpInf\_EP}. The system dynamics are accurately predicted via the $L_2$-regularized solution of \eqref{seqls3}.}
  \label{fig:states}
\end{subfigure}%
\hfill
\begin{subfigure}[t]{.48\textwidth}
  \centering
  \includegraphics[width=1.1\linewidth]{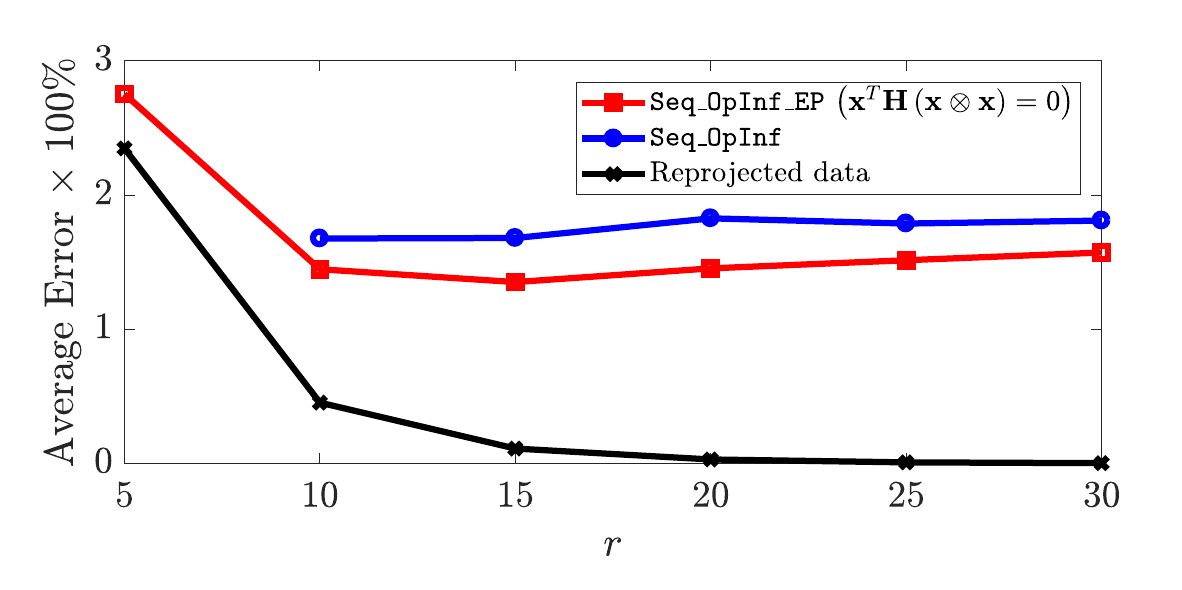}
    \caption{Average (in space and time) flowfield error for \eqref{burg}, for \texttt{Seq\_OpInf} and \texttt{Seq\_OpInf\_EP}. The \texttt{Seq\_OpInf} model for $r=5$ is unstable for the employed $L_2$ regularization. As a reference, the average error for the reprojected flowfield, $\mathbf{V} \mathbf{V}^\top \mathbf{x}$ is given.}
    \label{fig:error}
  \end{subfigure}
\vfill
\begin{subfigure}{1.0\textwidth}
  \centering
  \includegraphics[width=0.53\linewidth]{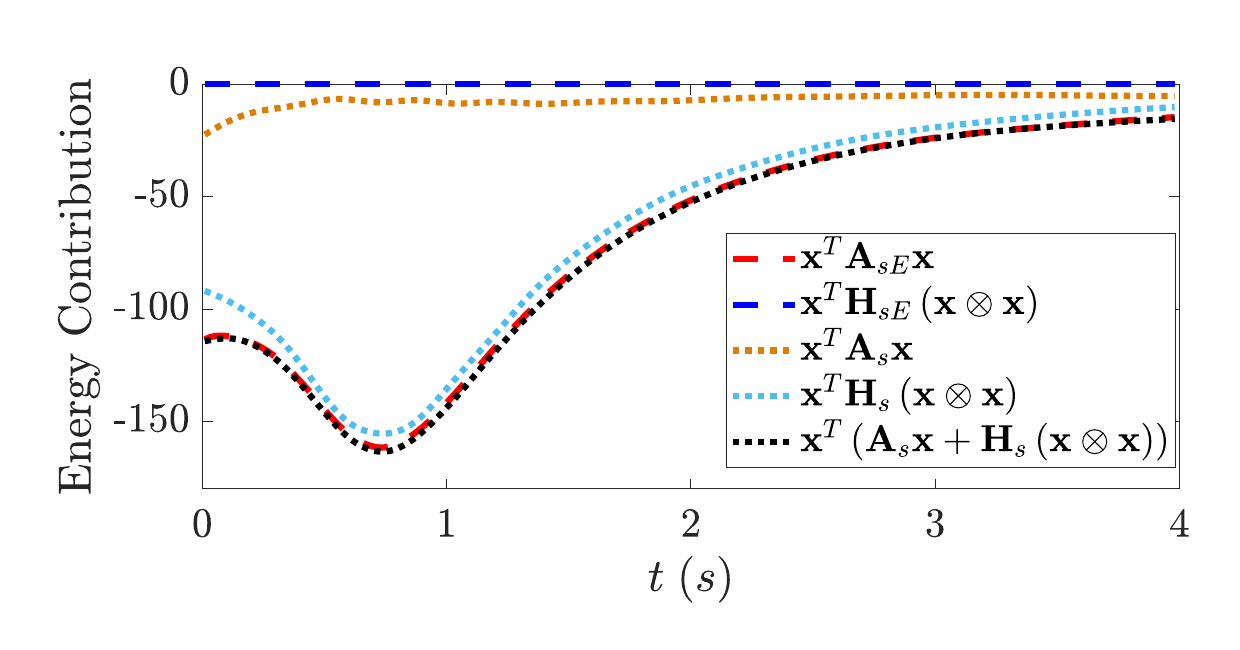}
    \caption{Energy norm over time, for the linear and quadratic term of \texttt{Seq\_OpInf} $(\mathbf{A}_s, \; \mathbf{H}_s)$ and \texttt{Seq\_OpInf\_EP} $(\mathbf{A}_{sE}, \; \mathbf{H}_{sE})$: Though in total, energy dissipation is correctly predicted by both models, standard \texttt{Seq\_OpInf} assigns energy dissipation principally to the quadratic term, departing from the physics-based property of $\mathbf{x}^\top\mathbf{H}_a \left( \mathbf{x} \otimes \mathbf{x} \right) = 0$. }
    \label{fig:energy}
  \end{subfigure}
\label{fig:test}
\caption{Numerical results for 2D Burgers' equation: Enforcing the physics-based energy preservation for the quadratic term through \texttt{Seq\_OpInf\_EP} produces accurate and more robust models in comparison to standard \texttt{Seq\_OpInf}.}
\end{figure}

The state prediction by the proposed \texttt{Seq\_OpInf\_EP} approach at $t=4$, with $r=15$ is shown in \Cref{fig:states}. The predictions in \Cref{fig:states} validate the accuracy of the method. The prediction error, normalized by the maximum state value, is also computed over time. The average error over space and time for \texttt{Seq\_OpInf} and \texttt{Seq\_OpInf\_EP} is given in \Cref{fig:error}, for different ROM dimensions. Exploiting the physics-based knowledge for the energy-preserving nonlinearity is beneficial, with \texttt{Seq\_OpInf\_EP} slightly out-performing standard OpInf in \texttt{Seq\_OpInf}.
 
Even though the predictions by both \texttt{Seq\_OpInf\_EP} and \texttt{Seq\_OpInf} are accurate, there is a qualitative difference between the two inferred models. This is indicated by \Cref{fig:energy}, where the energy contribution over time for the inferred linear term ($\mathbf{A}$) and the quadratic term ($\mathbf{H}$), for both inferred ROMs with $r=15$, is given. We denote the inferred operators for \texttt{Seq\_OpInf} by $\mathbf{A}_s$, $\mathbf{H}_s$ and the inferred operators for \texttt{Seq\_OpInf\_EP} by $\mathbf{A}_{sE}$, $\mathbf{H}_{sE}$. The cumulative energy dissipation by both the linear $\mathbf{A}_s$ and the quadratic $\mathbf{H}_s$ term (for \texttt{Seq\_OpInf}) is accurate, and essentially overlaps with the energy dissipation for $\mathbf{A}_{sE}$ and $\mathbf{H}_{sE}$ (for \texttt{Seq\_OpInf\_EP}). This is expected, since accurate predictions are made by both approaches for $r=15$ (see \Cref{fig:error}). However, the surprising observation is that the inferred ROM by \texttt{Seq\_OpInf} attributes energy dissipation mainly to the quadratic term ($\mathbf{H}_s$), while the contribution of the linear term ($\mathbf{A}_s$) to the system kinetic energy is minor. Standard OpInf, without the physics-based knowledge of energy-preserving nonlinearity, does neither correctly infer the dissipative properties expected by the linear term $\mathbf{A}_v$, or the energy-preservation behavior of $\mathbf{H}_a$ in \eqref{burg2}. On the contrary, exploiting the physics-informed property of energy-preserving nonlinearity due to \eqref{integral}, via the admissible representation \eqref{skew1} proven in \Cref{sec:fact}, the inferred ROM in \texttt{Seq\_OpInf\_EP} identically satisfies $\mathbf{x}^\top\mathbf{H}_{sE} \left( \mathbf{x} \otimes \mathbf{x} \right) = 0$. These additional, physics-informed constraints to \eqref{seqls5} allow to also capture the correct dissipative behaviour by $\mathbf{A}_{sE}$. We thus conclude that the inferred ROM by \texttt{Seq\_OpInf\_EP} is more robust and faithful to the physics of the underlying system.

\noindent The data and codes for this numerical example are available in \url{https://github.com/lgkimisis/quad_energy_opinf}.

\section{Discussion}
The main aim of this work was to prove that any energy-preserving quadratic operator can be represented by an equivalent operator with skew-symmetric sub-matrices. We showcased how this result can be employed in non-intrusive ROMs via the Operator Inference framework and presented a solution approach through sequential, linear LS problems. Results on a benchmark 2D Burgers' problem indicate that this approach produces accurate results, faithful to the underlying physics, with essentially no additional computational cost. An open research question considers the strategy for sequentially solving these problems; we are interested in examining the effect of the $j$-th LS solution to the subsequent, $(j+1)$-th LS problem, considering the inference of the remaining quadratic entries. Finally, the extent to which the proven representation can be employed to the inference of systems with only approximately energy-preserving nonlinearities lies in the future scope of this work.

\addcontentsline{toc}{section}{References}
\bibliographystyle{plainurl}
\bibliography{exampleref}

\appendix 
\section{Alternative proof}\label{appendix}
\begin{theorem*}\label{thm:equvilence_between_H}
	Consider an energy-preserving matrix $\bH \in \R^{n\times n^2}$, satisfying \eqref{energy_pres}. Then, there exists a matrix $\tilde \bH$ with structure as follows:
	\begin{equation}\label{eq:skew_symmetric_H}
	\tilde \bH = [\bH_1,\ldots, \bH_n],\qquad \bH_i = -\bH_i^\top  ~~\text{for}~~i = \{1,\ldots, n\},
	\end{equation}
	which is also energy-preserving and $\bH(\bx \otimes \bx) = \tilde \bH(\bx\otimes \bx)$ for every $\bx\in \R^n$.
\end{theorem*}

\begin{proof}
In this proof, we will show that a generic $\tilde \bH$ as in \eqref{eq:skew_symmetric_H} for which  $\bH(\bx \otimes \bx) = \tilde \bH(\bx\otimes \bx)$ is encoded by a linear system that has at least one solution. First note that $\tilde \bH$ satisfying \eqref{skew} is energy-preserving since $\bx^\top \tilde \bH (\bx\otimes \bx) = 0$. 

Since both matrices $\bH$ and $\tilde \bH$ are energy-preserving, from  \eqref{energy_pres}, we have
\begin{subequations}\label{eq:energy_preserving_appendix}
	\begin{align}
		e_k^\top\bH(e_i\otimes e_j + e_j\otimes e_i) + e_j^\top\bH(e_i\otimes e_k + e_k\otimes e_i) + 	e_i^\top\bH(e_j\otimes e_k + e_k\otimes e_j) &= 0, \\		e_k^\top\tilde\bH(e_i\otimes e_j + e_j\otimes e_i) + e_j^\top\tilde\bH(e_i\otimes e_k + e_k\otimes e_i) + 	e_i^\top\tilde\bH(e_j\otimes e_k + e_k\otimes e_j) &= 0,
	\end{align}
\end{subequations}
for  $i,j,k \in \{1,\ldots,n\}$.
Next, we notice that 
\begin{equation}
	\bH(\bx \otimes \bx) = \tilde \bH(\bx\otimes \bx)
	\end{equation}

 	\noindent if and only if
	\begin{equation}\label{eq:equvi_H}
	e_k^\top \bH(e_i \otimes e_j + e_j \otimes e_i) = e_k^\top \tilde\bH(e_i \otimes e_j + e_j \otimes e_i), ~~\text{for}~~{(i,j,k)} \in \{1,\ldots,n\}.
	\end{equation}

	Now, consider the case when $(i,j,k) \in \cM_1$, where $\cM_1 = \{(i, j, k) \mid i, j, k \in \{1, 2, \ldots, n\} \text{ and } i =j=k\}$. For this case, since $\tilde\bH$ and $\bH$ are energy-preserving, we have
	$$e_k^\top\bH(e_i \otimes e_j + e_j \otimes e_i) = e_k^\top\tilde\bH(e_i \otimes e_j + e_j \otimes e_i) = 0\qquad (i,j,k)\in \cM_1.$$
	
	Next, we focus on the cases when $(i,j,k) \in \cM_2$, where $\cM_2 = \{(i, j, k) \mid i, j, k \in \{1, 2, \ldots, n\} \text{ and } i \neq j \text{ or } i\neq k \text{ or } j\neq k\}$. First note that
	\begin{equation}
	e_k^\top \bH(e_i \otimes e_j + e_j \otimes e_i) = \left((e_i^\top \otimes e_j^\top + e_j^\top \otimes e_i^\top) \otimes e_k^\top\right)\bH_{\texttt{vect}},
	\end{equation}
	where $\bH_{\texttt{vect}}$ is the vectorization of the matrix $\bH$. Then, the conditions in \eqref{eq:equvi_H} for $(i,j,k) \in \cM_2$ can be written as follows:
	\begin{equation}\label{eq:equvi_H_Aform}
	\cA \bH_{\texttt{vect}} = \cA \tilde\bH_{\texttt{vect}},
	\end{equation}
	where $\cA$ is a linear operator defined so that its row vectors are constructed using $\left((e_i^\top \otimes e_j^\top + e_j^\top \otimes e_i^\top) \otimes e_k^\top\right)$ for all $(i,j,k) \in \cM_2$. Consequently, the matrix $\cA$ takes the following form:
	\begin{equation}
	\cA = \begin{bmatrix}
	\left((e_1^\top \otimes e_1^\top + e_1^\top \otimes e_1^\top) \otimes e_2^\top\right)\\
	\left((e_1^\top \otimes e_2^\top + e_2^\top \otimes e_1^\top) \otimes e_2^\top\right)\\
	\vdots \\
	\left((e_n^\top \otimes e_n^\top + e_n^\top \otimes e_n^\top) \otimes e_2^\top\right)\\
	\left((e_1^\top \otimes e_1^\top + e_1^\top \otimes e_1^\top) \otimes e_3^\top\right) \\
	\vdots \\ 
	\left((e_n^\top \otimes e_n^\top + e_n^\top \otimes e_n^\top) \otimes e_3^\top\right) \\
	\vdots \\
	\left((e_n^\top \otimes e_{n-1}^\top + e_{n-1}^\top \otimes e_n^\top) \otimes e_n^\top\right).
	\end{bmatrix}
	\end{equation}
	We note that the total number of rows in the matrix $\cA$  is $\tfrac{n^2(n+1)}{2} - n$. 
	
	Furthermore, we write the energy-preserving conditions \eqref{eq:energy_preserving_appendix} in a vectorized form for $\{i,j,k\} \in \cM_2$, thus yielding
	\begin{equation}\label{eq:SystemEP}
	\begin{aligned}
	\left((e_i\otimes e_j + e_j\otimes e_i)^\top \otimes e_k^\top  +  (e_i\otimes e_k + e_k\otimes e_i)^\top \otimes e_j^\top  + (e_j\otimes e_k + e_k\otimes e_j)^\top \otimes e_i^\top  \right)\bH_{\texttt{vect}}&= 0.
	\end{aligned}
	\end{equation}
    Equations as \eqref{eq:SystemEP} can be written in matrix form as:
	\begin{equation}
	\cB\bH_{\texttt{vect}} = 0,
	\end{equation}
	where $\cB$ can be constructed as follows:
	\begin{equation}
	\cB = \begin{bmatrix}
	\left((e_1\otimes e_1 + e_1\otimes e_1)^\top \otimes e_2^\top +  (e_1\otimes e_2 + e_2\otimes e_1)^\top \otimes e_1^\top + (e_1\otimes e_2 + e_2\otimes e_1)^\top \otimes e_1^\top \right) \\ 
	\left((e_2\otimes e_1 + e_2\otimes e_1)^\top \otimes e_2^\top +  (e_2\otimes e_2 + e_2\otimes e_2)^\top \otimes e_1^\top + (e_2\otimes e_1 + e_1\otimes e_2)^\top \otimes e_2^\top \right) \\ 
	\vdots \\
	\left((e_1\otimes e_2 + e_2\otimes e_1)^\top \otimes e_3^\top +  (e_1\otimes e_3 + e_3\otimes e_1)^\top \otimes e_2^\top + (e_2\otimes e_3 + e_3\otimes e_2)^\top \otimes e_1^\top \right) \\ 	
	\end{bmatrix}.
	\end{equation}
	Note that the total number of rows in the matrix $\cB$  is $n(n-1) + \tfrac{n(n-1)(n-2))}{6}$. Next, we find that the rows of the matrix $\cB$ are linearly independent and can be written as a linear combination of rows of the matrix $\cA$. Consequently, we can write $\bP\cA = \begin{bmatrix} \cA_1 \\ \cB\end{bmatrix}$, where $\bP$ is an invertible matrix and $\cA_1$ contains a subset of the rows of the matrix $\cA$. Then, after multiplying \eqref{eq:equvi_H_Aform} by the matrix $\bP$ on both sides, we get
	\begin{equation}
	\begin{aligned}
	\bP \cA \bH_{\texttt{vect}} &= 	\bP \cA \tilde\bH_{\texttt{vect}},\\
	\begin{bmatrix}
	\cA_1 \\ \cB
	\end{bmatrix} \bH_{\texttt{vect}} &= \begin{bmatrix}	\cA_1 \\ \cB \end{bmatrix} \tilde\bH_{\texttt{vect}}.
	\end{aligned}
	\end{equation}
	Due to energy-preserving properties of the matrices $\bH$ and $\tilde \bH$, we have $\cB \bH_{\texttt{vect}} =  \cB\tilde \bH_{\texttt{vect}} = 0$. Hence, to ensure $\bH(\bx\otimes \bx) = \tilde{\bH}(\bx\otimes \bx)$ for every $\bx \in \R^n$, provided the matrices $\bH$ and $\tilde \bH$ are energy-preserving, we only require the following condition to be satisfied:
	\begin{equation}\label{eq:constraint_1}
	\begin{aligned}
	\cA_1 \bH_{\texttt{vect}} &= \cA_1  \tilde\bH_{\texttt{vect}},
	\end{aligned}
	\end{equation}
	with the number of rows of the matrix $\cA_1$ to be $$\dfrac{n^2(n+1)}{2} - n - n(n-1) - \dfrac{n(n-1)(n-2))}{6} = \dfrac{n(n+1)(n-1)}{3}.$$

Additionally we impose an additional structure on the matrix $\tilde \bH$, particularly the one given in \eqref{eq:skew_symmetric_H}. We begin by noticing that due to the particular skew-symmetric structure, we have $\dfrac{n^2(n+1)}{2}$ constraints on the matrix $\tilde \bH$, which can be written as follows:
	\begin{equation}
	\left(e_k^\top\otimes (e_i\otimes e_j + e_j\otimes e_i)\right)\tilde{\bH}_{\texttt{vect}}  = 0,
	\end{equation}
	where $(i,j,k) \in \{1,\ldots,n\}$. We can write these constraints in matrix form as follows:
	\begin{equation}
	\cC \tilde{\bH}_{\texttt{vect}} = 0.
	\end{equation} 
	Combining the above equation with \eqref{eq:constraint_1}, we obtain 
	\begin{equation}\label{eq:solve_H_tH}
	\begin{bmatrix}
	\cA_1 \\ \cC
	\end{bmatrix} \tilde{\bH}_{\texttt{vect}} = \begin{bmatrix} 	
	\cA_1 \tilde{\bH}_{\texttt{vect}} \\ 0 \end{bmatrix}.
	\end{equation}
	Next, note that the rows of the matrices $\cA_1$ and $\cC$ are linearly independent. The reason is that the rows are $\cA_1$ are constructed using 
	\begin{equation}
	\left((e_i^\top \otimes e_j^\top + e_j^\top \otimes e_i^\top) \otimes e_k^\top\right)\qquad {(i,j,k)} \in \cM_2
	\end{equation}
	and the rows are $\cC$ are constructed using 
	\begin{equation}
	e_k^\top\otimes \left((e_i^\top \otimes e_j^\top + e_j^\top \otimes e_i^\top)\right)\qquad {(i,j,k)} \in \cM_1 \cup \cM_2,
	\end{equation}
	and they are linearly independent. This follows from the property that $e_k^\top \otimes \ba$ and $\ba \otimes e_k^\top $ are linearly independent for every $\ba$ and $k \in \{1,\ldots, n\}$ provided $\ba\neq \alpha e_k$ for $\alpha \in \R$, and the vectors corresponding to $\ba = \alpha e_k$ are not present in the rows of $\cA_1$. Therefore, the matrix $		\begin{bmatrix}
	\cA_1 \\ \cC
	\end{bmatrix}$ is of full-row rank with $$\dfrac{n^2(n+1)}{2} + \dfrac{n(n+1)(n-1)}{3}$$ rows, less than (for $n>2$) or equal (for $n\in\{1 , 2\}$) to $n^3$. Hence, a solution always exists to the equation \eqref{eq:solve_H_tH}. This ensures the existence of $\tilde \bH$ in the form \eqref{eq:skew_symmetric_H} such that $\bH(\bx\otimes \bx) = \tilde \bH(\bx\otimes \bx)$ for all $\bx \in \Rn$ provided that $\bH$ is energy-preserving. This concludes the proof.
\end{proof}  
\end{document}